\documentclass[12pt]{article}
\usepackage{amsmath,amssymb,float,tikz}
\usepackage{graphicx}

\begin{document}


\begin{center}
{\large\bf Santalo's formula and stability of trapping sets\\ of positive measure}
\end{center}

\begin{center}
{\bf Luchezar Stoyanov\footnote{E-mail: luchezar.stoyanov@uwa.edu.au}}
\end{center}


\begin{center}
{\footnotesize\it School of Mathematics and Statistics, University of Western Australia, Crawley 6009 WA, Australia}
\end{center}

\medskip

\newcommand{\R}{{\sf I\hspace{-.15em}R}}
\def\C{{\bf C}}
\def\e{\emptyset}
\def\ds{\partial S}
\def\cq{\overline{Q}}
\def\sn{{\bf S}^{n-1}}
\def\hg{\Gamma}
\def\ssn{{\sn}\times {\sn}}
\def\do{\partial \Omega}
\def\dk{\partial K}
\def\dl{\partial L}
\def\tts{\tilde{S}}
\def\ttp{\tilde{p}}
\def\tu{\tilde{U}}
\def\hs{\hat{S}}
\def\hp{\hat{p}}
\def\dr{\frac{\partial r}{\partial x_1}}
\def\ll{{\cal L}}
\def\hm{\hat{M}}
\def\tM{\widetilde{M}}
\def\pp{{\cal P}}
\def\tt{{\cal T}}
\def\ds{\partial S}
\def\ss{{\cal S}}
\def\vv{{\cal V}}
\def\aa{{\cal A}}
\def\bb{{\cal B}}
\def\nn{{\cal N}}
\def\dd{{\cal D}}
\def\ee{{\cal E}}
\def\uu{{\cal U}}
\def\rr{{\cal R}}
\def\kk{{\cal K}}
\def\cc{C_0^{\infty}}
\def\ot{(\omega,\theta)}
\def\oto{(\omega_0,\theta_0)}
\def\toto{(\tilde{\omega}_0, \tilde{\theta}_0)}
\def\ts{\tilde{\sigma}_0}
\def\tx{\tilde{x}}
\def\txi{\tilde{\xi}}
\def\ttt{\tilde{t}}
\def\got{\gamma(\omega,\theta)}
\def\ggot{\gamma'(\omega,\theta)}
\def\omo{(\omega,\omega)}
\def\du{\partial u_i}
\def\dfu{\frac{\partial f}{\partial u_i}}
\def\dou{\frac{\partial \omega}{\partial u_i}}
\def\dttu{\frac{\partial \theta}{\partial u_i}}
\def\ooo{(\omega_0,-\omega_0)}
\def\otp{(\omega',\theta')}
\def\oo{{\cal O}}
\def\pr{{\rm pr}}
\def\cK{\hat{K}}
\def\cL{\hat{L}}
\def\ff{{\cal F}}
\def\fk{{\cal F}^{(K)}}
\def\fl{{\cal F}^{(L)}}
\def\kkr{{\cal K}^{\mbox{reg}}}
\def\kkro{{\cal K}_0^{\mbox{reg}}}
\def\trapk{\mbox{\rm Trap}(\Omega_K)}
\def\trap{\mbox{\rm Trap}}
\def\trapo{\mbox{\rm Trap}_\Omega}
\def\trapf{\mbox{\footnotesize\rm Trap}}
\def\hr{\hat{\rho}}
\def\G{{\cal G}}
\def\tg{\tilde{\gamma}}
\def\Vol{\mbox{Vol}}
\def\dkso{\partial K^{(\infty)}}
\def\dlso{\partial L^{(\infty)}}
\def\gg{{\cal g}}
\def\sl{{\cal SL}}
\def\i{{\bf i}}
\def\te{{\cal T}^{(ext)}}
\def\endofproof{{\rule{6pt}{6pt}}}
\def\Box{\endofproof}
\def\su{S^*(\R^n\setminus U)}
\def\gk{\gamma_K}
\def\gl{\gamma_L}
\def\la{\left\langle}
\def\ra{\right\rangle}
\def\kfin{\kk_0^{({\mbox{\footnotesize\rm fin}})}} 
\def\dt{\dot{T}}
\def\ep{\epsilon}
\def\kfi{\kk^{({\mbox{\footnotesize\rm fin}})}} 
\def\stk{\Sigma_3^{(K)}}
\def\stl{\Sigma_3^{(L)}}
\def\SU{S^*(\R^n\setminus U)}
\def\tkm{\tt_k^{(m)}}
\def\ukm{U_k^{(m)}}
\def\Pkm{\Psi_k^{(m)}}
\def\nkm{N_k^{(m)}}
\def\di{\displaystyle}
\def\ep{\epsilon}
\def\dist{\mbox{\rm dist}}
\def\diam{\mbox{\rm diam}}
\def\con{\mbox{\rm const}}
\def\bs{\bigskip}
\def\ms{\medskip}
\def\id{\mbox{\rm id}}
\def\SK{S^*_{K}(S_0)}
\def\SL{S^*_{L}(S_0)}
\def\Ir{I^{(r)}}
\def\dM{\partial M}

\def\te{\tilde{e}}
\def\tx{\tilde{x}}
\def\ty{\tilde{y}}
\def\grad{\; \mbox{\rm grad} \,}
\def\spp{\vspace{5pt}\noindent}
\newtheorem{definition}{Definition}
\newtheorem{construction}{Construction}
\newtheorem{example}{Example}
\newtheorem{lemma}{Lemma}
\newtheorem{theorem}{Theorem}
\newtheorem{corollary}{Corollary}
\newtheorem{proposition}{Proposition}

\footnotesize

\noindent
{\bf Abstract} 

Billiard trajectories (broken generalised geodesics) are considered in the exterior of an obstacle 
$K$ with smooth boundary on an arbitrary Riemannian manifold. We prove a generalisation of the well-known
Santalo's formula. As a consequence, it is established that if the set of trapped points
has positive measure, then for all sufficiently small smooth perturbations of the boundary of $K$ the set of trapped 
points for the new obstacle $K'$ obtained in this way also has positive measure. More generally 
the measure of the set of trapped points depends continuously on perturbations of the obstacle $K$. Some
consequences are derived in the case of scattering by an obstacle $K$ in $\R^n$. For example, it is shown that, 
for a large class of obstacles $K$, the volume of $K$ is uniquely determined by the average travelling times
of scattering rays in the exterior of $K$.

\medskip 


\medskip


\normalsize

\section{Introduction}\label{sec1}
\renewcommand{\theequation}{\arabic{section}.\arabic{equation}}

Let $g_t$ be the {\it geodesic flow} on the {\it unit sphere bundle} $S\tM$ of a $C^k$ ($k \geq 2$) Riemannian manifold $\tM$
without boundary, and let $n = \dim(\tM) = n \geq 2$. Consider an arbitrary compact subset $M$ of $\tM$ with non-empty smooth 
($C^k$) boundary $\partial M$ and non-empty interior $M\setminus \dM$. 

Next, suppose that $K$ is a compact subset of the interior 
$M \setminus \partial M$ of $M$ with smooth boundary $\dk$ and non-empty interior $K \setminus \dk$. 
Consider the compact subset
$$\Omega = \overline{M\setminus K}$$
of $M$. Then the boundary of $\Omega$ in $\tM$ is $\do = \partial M \cup \dk$. Define the {\it billiard flow}
$\phi_t$ in $S(\Omega)$ in the usual way:  it coincides with the geodesic flow $g_t$ in the interior of $S(\Omega)$, and when a geodesic hits the 
boundary $\do$ it reflect following the law of geometrical optics. Let $dq$ be the measure on $\tM$ determined by the Riemannian metric and let
$dv$ be the Lebesgue measure on the unit sphere $S_q\tM$. It is well-known (see Sect. 2.4 in \cite{CFS})
that the geodesic flow $g_t$ preserves the Liouville  measure $dq dv$ on $S\tM$, and so the billiard flow $\phi_t$ preserves the {\it restriction} 
$\lambda$ of $dq dv$ to $S(\Omega)$. Finally, denote by $\mu$ the {\it Liouville measure } on $S^+(\do)$ defined by 
$$d\mu = d\rho(q) d\omega_q |\la \nu(q), v\ra| ,$$
where $\rho$  is the measure on $\do$ determined by the Riemannian metric on $\do$ and $\omega_q$ is the Lebesgue measure on the 
$(n-1)$-dimensional unit sphere $S_q(M)$ (see e.g. Sect. 6.1 in \cite{CFS}).

For any $q\in \do$, denote by $\nu(q) \in S\tM$ the {\it unit normal vector} to $\do$ pointing into the interior of $\Omega$. Set
$$S^+(\do) = \{ x = (q,v) : q \in \do , \la v, \nu(q) \ra \geq 0\} ,$$
and in a similar way define $S^+(\dM) \subset S^+(\do)$. For any $x = S^+(\do)$ denote by $\tau(x) \geq 0$ the maximal number (or $\infty$)
such that $\phi_t(x) = g_t(x)$ is in the interior of $\Omega$ for all $0 < t < \tau(x)$. Set $\tau(x) = 0$ if
$\la \nu(q),v\ra = 0$. The function $\tau$ thus defined is called the {\it first return time} function for the flow $\phi_t$ on $\Omega$.

In the present paper we will assume that the geodesic flow $g_t$ has no trapped trajectories  in $M$, i.e.

\ms

\noindent
{\bf (A):} {\it For every $x = (q,v) \in S^+(\dM)$ with $\la v, \nu(q)\ra > 0$ $\exists t > 0$ with $\pr_1(g_t(x)) \in \dM$,}

\ms

\noindent
where we set $\pr_1(q,v) = q$. In particular, it follows from this assumption that $\tau(x) < \infty$ for every $x\in S^+(\do)$. 

Then {\it Santalo's formula} (see e.g. pp. 336-338 in \cite{Sa} or \cite{Cha}) gives
\begin{equation}
\int_{S(\Omega)} f(x) \; d\lambda(x) = \int_{S^+(\do)} \left(\int_0^{\tau(x)} f(g_t(x)) \;dt \right)\, d\mu(x)
\end{equation}
for every $\lambda$-integrable function $f$ on $S(\Omega)$.

Santalo's formula has been widely used in geometry and dynamical systems (including billiards) -- 
see e.g. \cite{CULV}, \cite{Ch}, \cite{CaG}, \cite{C} and the references there. 

By assumption (A), the geodesic flow $g_t$ has no trapped trajectories in $M$.
However, the billiard flow $\phi_t$ may have trapped trajectories in $M$ with respect to the {\it obstacle} $K$.
More precisely, we will now consider the billiard trajectories in $\Omega$ as scattering 
trajectories in $M$ reflecting at the boundary $\dk$  when they hit $K$. Given $x = S(\Omega)$, let
$t(x) \geq 0$ be the maximal number (or $\infty$)  such that $\phi_t(x)$ is in the 
interior of $M$ for all $0 < t < t(x)$. For $x = (q,v) \in S^+(\do)$ with $\la \nu(q),v\ra = 0$ set $t(x) = 0$. 
Then $t(x)$ is called the {\it travelling time} function on $\Omega$. Given $x\in S^+(\partial M)$, consider
the billiard trajectory 
$$\gamma^+(x) = \{ \phi_t(x) : 0 \leq t \leq t(x)\} $$
which starts at a point  $x\in S^+(\partial M)$ and, if $t(x) < \infty$, it ends at another point $y \in S^+(\partial M)$. 
Between $x$ and $y$ the trajectory may hit $\dk$ several times. 

Problems related to recovery of information about the manifold $M$ from travelling times $\tau(x)$ of geodesics in the interior of $M$
have been considered in Riemannian geometry for a rather long time -- see \cite{SUV}, \cite{SU}, \cite{CULV} and the references there for some 
general information. Similar problems have been studied when an obstacle $K$ is present; then one deals with travelling times $t(x)$ of 
billiard  trajectories (generalised geodesics) in the exterior of $K$ -- see  \cite{NS1}, \cite{NS2},
\cite{St3} and the references there (see also \cite{LP} and \cite{M} for some general information on 
inverse scattering problems).  For  some classes of obstacles $K$ the travelling time function $t(x) = t_K(x)$ completely determines $K$ 
(\cite{NS1}, \cite{NS2}). For example, it was recently established (\cite{NS1}) that this is so in the class 
of obstacles $K$ in Euclidean spaces that are finite disjoint unions of strictly convex domains with smooth boundaries.

Let $\trapo^+(\dM)$ be the set of the {\it trapped points} $x = (q,v) \in S^+(\dM)$ of the billiard flow $\phi_t$ in $\Omega$, 
i.e. the points for which $t(x) = \infty$. In general, it may happen that $\trapo^+(\dM) \neq \e$, however it follows from 
Theorem 1.6.2 in \cite{LP} (see also Proposition 2.4 in \cite{St1}) that 
\begin{equation}
\mu(\trapo^+(\dM)) = 0 .
\end{equation}

Next, let $\trap(\Omega)$ be the set of all {\it trapped points} $x = (q,v) \in S(\Omega)$ of the billiard flow $\phi_t$ in $\Omega$, 
i.e. the points for which  $t (q,v) = \infty$ or $t (q,-v) = \infty$. As Livschits' example shows (see Figure 1), {\bf in general $\trap(\Omega)$ 
may have positive  Liouville measure and even a non-empty interior}.

Generalising  Santalo's formula (1.1), we prove the following.

\bs

\noindent
{\bf Theorem 1.1.} { \it Assume that $M$ and the geodesic flow $g_t$ in $M$ satisfy the condition {\rm (A)}.  
Let  $f : S(\Omega) \setminus \trap(\Omega)\longrightarrow \C$ be a $\lambda$-measurable function.
Assume that either:}

 (i) {\it $\trap(\Omega) = \e$  and $f$ is integrable on $S(\Omega)$, or }
 
 (ii) {\it $|f|$ is integrable.}
 
 \noindent
{\it Then  we have}
\begin{equation}
\int_{S(\Omega)\setminus \trapf(\Omega)} f(x) \; d\lambda(x) 
= \int_{S^+(\partial M) \setminus \trapf^+_{\Omega}(\dM)} \left( \int_0^{t(x)} f(\phi_t(x)) \;dt \right)\, d\mu(x) .
\end{equation}

\ms

Clearly, in the case (ii)  above $\trap(\Omega)$ is allowed to  have positive measure.   

We will now describe the main consequence of Theorem 1.1 derived in this paper.

Let $k \geq 3$ and let $C^k(\dk, M)$ be {\it the space of all smooth embedding} $F : \dk \longrightarrow M$ endowed with the Whitney
$C^k$ topology (see \cite{Hi}). Given $F \in C^k(\dk, M)$, let $K_F$ be the obstacle in $M$ with boundary $\dk_F = F(\dk)$ so that
$K_F \cap \dM = \e$, and let $\Omega_F = \overline{M\setminus K_F}$.

Our second main result in this paper is the following.

\bs

\noindent
{\bf Theorem 1.2.} { \it Assume that $M$ and the geodesic flow $g_t$ in $M$ satisfy the condition {\rm (A)}.  }

\ms

(a) {\it If $\lambda(\trap(\Omega)) > 0$, then there exists an open neighbourhood $U$ of $\id$ in $C^k(\dk,M)$ such that for 
every $F \in U$ we have $\lambda (\trap(\Omega_F)) > 0$.}

\ms

(b) {\it More generally, for every $\ep > 0$  there exists an open neighbourhood $U$ of $\id$ in $C^k(\dk,M)$ such that for 
every $F \in U$ we have $|\lambda (\trap(\Omega_F)) - \lambda(\trap(\Omega))| < \ep$.}

\bs

We prove Theorem 1.1 in Sect. 2 below, and then use it in Sect.3 to derive  Theorem 1.2.
Let us note that in the case $\trap(\Omega) = \e$ the formula (1.3) with $f = 1$ was mentioned without proof and used by 
Plakhov and Roshchina in \cite{PlaR} (see also \cite{Pla}).

Formula (1.3)  with $f = 1$ implies the following.

\bs

\noindent
{\bf Corollary 1.3.} { \it Under the assumptions of Theorem {\rm 1.1} we have}
$$\lambda (\trap(\Omega)) = \lambda(S(\Omega)) -  \int_{S^+(\partial M)} t(x)\, d\mu(x) .$$

\bs

So, if we know the travelling time function $t(x)$ and have enough information about $\Omega$ to determine its volume,
then we can determine the measure of the set of trapped points in $\Omega$, as well.

Some other consequences of Theorem 1.1 concerning scattering by obstacles are discussed in Sects. 4 and 5  below.

\bs

\begin{tikzpicture}[xscale=1.70,yscale=0.90] 

  \draw[violet,thick] (-1.5,3) arc (0:180:2cm and 2cm); 
  \draw (0,1) arc (-30:210:4cm and 4cm); 
  \draw (-3.8,2.3) arc (0:180:0.5cm and 0.7cm); 
  \draw (-2.1,2.3) arc (0:180:0.5cm and 0.7cm); 
  \draw (-3,5) node[anchor=south west] {$E$}; 
  \draw (-4.4,3) node[anchor=north west] {$F_1$};
  \draw (-2.4,3) node[anchor=north east] {$F_2$};
    \draw (-5.9,3.3) node[anchor=north west] {$P$};
  \draw (-1.0,3.3) node[anchor=north east] {$Q$ };
    \draw (-5.0,1.8) node[anchor=south east] {${\huge U}$ };
       \draw (-1.65,1.8) node[anchor=south east] {${\huge V}$ };


\draw (-6.92,1) .. controls (-5.6,-1) and (-3.7,-0.7) .. (-3.8,2.29);
\draw (-3.1,2.29) .. controls (-3.2,-1) and (-1,-0.7) .. (0,1);
\draw[blue, thick] (-2.1,2.29) .. controls (-2,1.2) and (-1.5,1.5) .. (-1.5,3);
\draw[blue, thick] (-5.5,3.1) .. controls (-5.5,1.1) and (-4.9,1.4) .. (-4.8,2.4);

\draw[thick] (-5.5,3) circle (0.02cm);
\draw[thick] (-4.3,3) circle (0.02cm);
\draw[thick] (-2.6,3) circle (0.02cm);
\draw[thick] (-1.5,3) circle (0.02cm);

\end{tikzpicture}

\medskip

\begin{center}
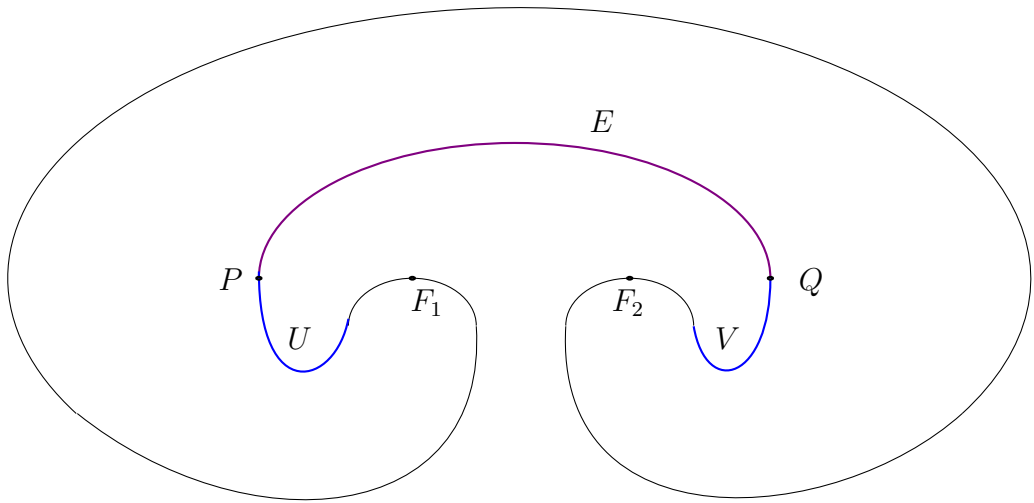
{\bf Figure 1: Livshits' Example} (Ch. 5 in \cite{M})
\end{center}
\noindent
{\footnotesize $E$ is a half-ellipse with end points $P$ and $Q$ and foci $F_1$ and $F_2$.
 Any trajectory entering the interior of the ellipse between the foci must exit between the foci after reflection. So,
no trajectory `coming from infinity' has a common point  with the parts $U$ and $V$ of the boundary and the nearby regions.}

\medskip

Similar examples in higher dimensions are described in \cite{NS3}.

\section{A generalised Santalo's formula}\label{sec2}
\renewcommand{\theequation}{\arabic{section}.\arabic{equation}}
\setcounter{equation}{0}

Given $x\in S(\Omega)$, we will say that
$\gamma^+(x)$ {\it contains  a gliding segment on the boundary $\do$} if there exist $0 \leq t_1 < t_2 \leq t(x)$ 
such that $\phi_t(x) = g_t(x) \in S (\do)$ for all $t \in [t_1,t_2]$ (i.e. $\gamma^+(x)$ contains a non-trivial 
geodesic segment lying entirely in $\do$).

It follows from results\footnote{See Sect. 3 in \cite{MS2}; see also \cite{MS1}, \cite{H} or \cite{PS} for general
information about generalised geodesics.} in \cite{MS2}  that 
for $\lambda$-almost all $x \in S(\Omega)$, the billiard trajectory $\gamma^+(x)$ does not contain any gliding segments on 
the boundary $\do$, and $\gamma^+(x)$ has only finitely many reflections.

\bs

First, we prove a special case of Theorem 1.1.

\bs

\noindent
{\bf Lemma 2.1.} {\it Assume that $M$ and the geodesic flow $g_t$ in $M$ satisfy the condition  {\rm (A)}. 
Let $V$ be an open subset of $S(\Omega)$ containing $\trap(\Omega)$ and such that $\phi_t(x) \in V$ for any $x\in V$ and any $t \in [0,t(x)]$. 
Assume that $f : S(\Omega) \longrightarrow \C$ is an integrable function such that  $f = 0$ on $V$. Then {\rm (1.3)} holds.}

\bs

\noindent
{\bf Proof of Lemma 2.1. }
Let $V$ and $f$ satisfy the assumptions of the lemma. Then we have
\begin{equation}
\int_{S(\Omega)} f(x) \; d\lambda(x) = \int_{S(\Omega)\setminus V} f(x) \; d\lambda(x) 
\end{equation}
and 
\begin{equation}
\int_{S^+(\partial M)} \left( \int_0^{t(x)} f(\phi_t(x)) \;dt \right)\, d\mu(x) 
= \int_{S^+(\partial M)\setminus V} \left( \int_0^{t(x)} f(\phi_t(x)) \;dt \right)\, d\mu(x) .
\end{equation}
We will prove that the right-hand-sides of (2.1) and (2.2) are equal.

Given an integer $k \geq 0$ denote by $\Gamma_k$ the set of those $x \in S^+(\partial M)\setminus V$  such that
$\gamma^+(x)$ contains no gliding segments on the boundary $\do$,  and $\gamma^+(x)$ has exactly 
$k$ reflections at $\dk$. Clearly $\Gamma_k$ are disjoint, measurable subsets of $S^+(\partial M)$ and, as remarked earlier,
it follows from \cite{MS2} that
\begin{equation}
\mu \left( S^+(\partial M) \setminus (V \cup \cup_{k=0}^\infty \Gamma_k )\right) = 0 .
\end{equation}

The {\it billiard ball map} $B$ is defined in the usual way: given $x = (q,v) \in S^+(\do)$ with $g_t(x) \in S(\do)$ for some
$t \in (0,\infty)$, take the smallest $t > 0$ with this property, and let $g_t(x) = (p,w)$ for some $p \in \do$,  $w\in \sn$. 
Set $B(x) = (p, \sigma_p(w)) \in S^+(\do)$, where 
$$\sigma_p : T_p(\tM) \longrightarrow T_p(\tM)$$ 
is the {\it symmetry} through  the tangent  plane $T_p(\do)$,  i.e.
$$\sigma_p(\xi) = \xi - 2\la \nu(p), \xi \ra \, \nu(p) .$$

It follows from the condition (A)  that $B$ is well-defined on
a set $\Lambda$ of full $\mu$-measure in $S^+(\do)$ and the Liouville measure $\mu$ is invariant with respect to $B$ 
(see e.g. Lemma 6.6.1 in \cite{CFS}).

Notice that each of the sets $\Gamma_k$ is contained in $\Lambda$ and moreover
$$\Gamma_{k,j} = B^j(\Gamma_k) \subset \Lambda$$ 
for all $j = 0,1, \ldots,k$. Also, by the definition of the sets $\Gamma_k$ we have
\begin{equation}
B^j(\Gamma_k) \subset  S^+(\dk) \quad, \quad 1 \leq j \leq k .
\end{equation}
The sets $\Gamma_{k,j}$ are clearly measurable, and moreover
\begin{equation}
\Gamma_{k,j} \cap \Gamma_{m,i} = \e \quad \mbox{\rm whenever}\: k \neq m \:\:\mbox{\rm or }\: j \neq i .
\end{equation}
Indeed, assume that $y\in \Gamma_{k,j} \cap \Gamma_{m,i}$ for some non-negative numbers $k,j,m,i$
with $0 \leq j \leq k$ and $0 \leq i \leq m$. 
Then $y =B^j(x)$ for some $x\in \Gamma_k$ and $y = B^i(z)$ for some $z \in \Gamma_m$.
Assume e.g. $j > i$. Then $B^i(z) = y = B^j(x) = B^i(B^{j-i}(x))$ implies $z = B^{j-i}(x)$. Now (2.4)
gives $z \in S^+(\dk)$ which is a contradiction with $z \in \Gamma_m \subset S^+(\partial M)$. Thus, 
we must have $j = i$ and therefore $B^i(z) = y = B^i(x)$, so $z = x \in \Gamma_k \cap \Gamma_m$.
The latter is non-empty only when $k = m$. This proves (2.5).

Finally, notice that
\begin{equation}
\mu \left( S^+(\do) \setminus (V \cup \cup_{k=0}^\infty \cup_{j=0}^k \Gamma_{k,j} ) \right) = 0 .
\end{equation}
Indeed, as mentioned earlier, it follows from results in \cite{MS2} that 
for $\lambda$-almost all $x \in S(\Omega)$, the billiard trajectory $\gamma^+(x)$ does not contain any gliding segments on 
the boundary $\do$, and $\gamma^+(x)$ has only finitely many reflections. Since $\trap(\Omega) \subset V$,
for almost all $y \in S^+(\dk) \setminus V$ there exist $t \geq 0$ and
$x \in S^+(\partial M)$ such that $y = \phi_t(x)$, $t(x) < \infty$, the billiard trajectory   
$\gamma^+(x)$ does not contain any gliding segments on $\do$ and has only finitely many reflections.
Let $k$ be the number of those reflections; then for some $j = 0,1, \ldots,k$ we have $y = B^j(x)$.
Thus, $x \in \Gamma_k$ and $y \in \Gamma_{k,j}$. This proves (2.6).

Next, given $k > 0$ and $x \in \Gamma_k$, clearly we have
$$t(x) = \tau(x) + \tau(B(x)) + \tau(B^2(x)) + \ldots + \tau(B^k(x)) .$$
Set 
$$T_j(x) = \tau(x) + \tau(B(x)) + \ldots + \tau(B^j(x))$$ 
for $j = 0, 1, \ldots,k$.  
Clearly $\phi_{T_j(x)}(x) = B^{j+1}(x)$ for $0 \leq j \leq k-1$.
Hence for any $j = 1, \ldots,k$ we have 
$$\phi_{s+T_{j-1}(x)} (x) = \phi_{s} (\phi_{T_{j-1}(x)}(x)) = g_s(B^j(x)) ,$$
and therefore, using the substitution $t = s+ T_{j-1}(x)$ below we get
$$\int_{T_{j-1}(x)}^{T_j(x)} f(\phi_t(x)) \;dt
= \int_{0}^{\tau(B^j(x))} f(\phi_s(B^j(x))) \;ds .$$

Thus,
\begin{eqnarray*}
\int_0^{t(x)} f(\phi_t(x)) \;dt 
& = & \int_0^{\tau(x)} f(\phi_t(x)) \;dt +
\sum_{j=1}^k \int_{T_{j-1}(x)}^{T_j(x)} f(\phi_t(x)) \;dt \\
& = & \int_0^{\tau(x)} f(\phi_t(x)) \;dt +
\sum_{j=1}^k \int_{0}^{\tau(B^j(x))} f(\phi_{s}(B^j(x))) \;ds \\
& = & \sum_{j=0}^k \int_{0}^{\tau(B^j(x))} f(\phi_{t}(B^j(x))) \;dt .
\end{eqnarray*}

Now, using the above and the fact that $B^j : \Gamma_{k,j} \longrightarrow B^j(\Gamma_{k,j})$ is a
measure preserving bijection, we get
$$\int_{\Gamma_k} \left( \int_{0}^{\tau(B^{j}(x))} f(\phi_t(B^j(x))) \;dt \right)\, d\mu(x)
= \int_{\Gamma_{k,j}} \left( \int_{0}^{\tau(y)} f(\phi_t(y)) \;dt \right)\, d\mu(y) ,$$
which implies
\begin{eqnarray*}
\int_{\Gamma_k} \left(\int_0^{t(x)} f(\phi_t(x)) \;dt \right)\, d\mu(x)
& = & \int_{\Gamma_k} \left( \sum_{j=0}^k \int_0^{\tau(B^j(x))} f(\phi_t(B^j(x))) \;dt \right)\, d\mu(x)\\
& = &  \sum_{j=0}^k \int_{\Gamma_k} \left( \int_{0}^{\tau(B^{j}(x))} 
f(\phi_t(B^j(x))) \;dt \right)\, d\mu(x)\\
& = &  \sum_{j=0}^k \int_{\Gamma_{k,j}} \left( \int_{0}^{\tau(y)} f(\phi_t(y)) \;dt \right)\, d\mu(y) .
\end{eqnarray*}

Using (2.2), (2.3), (2.4) and (2.6), we derive
\begin{eqnarray*}
\int_{S^+(\partial M)} \left( \int_0^{t(x)} f(\phi_t(x)) \;dt \right)\, d\mu(x)
& = & \int_{S^+(\partial M)\setminus V} \left( \int_0^{t(x)} f(\phi_t(x)) \;dt \right)\, d\mu(x)\\
& = & \sum_{k=0}^\infty \int_{\Gamma_k} \left( \int_0^{t(x)} f(\phi_t(x)) \;dt \right)\, d\mu(x)\\
& = & \sum_{k=0}^\infty  \sum_{j=0}^k\int_{\Gamma_{k,j}} \left( \int_{0}^{\tau(y)}  f(\phi_t(y)) \;dt \right)\, d\mu(y)\\
& = & \int_{S^+(\do)\setminus V} \left( \int_{0}^{\tau(y)} f(\phi_t(y)) \;dt \right)\, d\mu(y) \\
& = & \int_{S^+(\do)} \left( \int_{0}^{\tau(y)} f(g_t(y)) \;dt \right)\, d\mu(y) .
\end{eqnarray*}
By Santalo's formula (1.1), the latter is equal to $\di \int_{S(\Omega)} f(x) \; d\lambda(x)$. 
\endofproof

\bs

\noindent
{\bf Proof of Theorem 1.1.}  If $\trap(\Omega) = \e$,  taking $V = \e$ in Lemma 2.1, proves the case (i).

To prove the case (ii), assume that $\trap(\Omega) \neq \e$. 
Let $f$ be a measurable function of $S(\Omega) \setminus \trap(\Omega)$ such that $|f|$ is integrable.
Extend\footnote{This is already assumed in the right-hand-side of (1.3).} $f$ as $0$ on $\trap(\Omega)$. 
Since $\trap^+_\Omega(\dM)$ is closed in $S^+(\dM)$, there exists a sequence of  compact sets
$$F_1 \subset F_2 \subset \ldots \subset F_m \subset  \ldots \subset S^+(\dM)\setminus \trap^+_\Omega(\dM)$$
such that
\begin{equation}
\cup_{m=1}^\infty F_m = S^+(\dM)\setminus \trap^+_\Omega(\dM) .
\end{equation}
Choose such a sequence, and for every $m$ set
$$G_m = \{ \phi_t(x) : x\in F_m \:, \: 0\leq t \leq t(x) \} .$$
Then we have
\begin{equation}
\cup_{m=1}^\infty G_m = S(\Omega) \setminus \trap (\Omega) .
\end{equation}
Indeed, if $x \in G_m$ for some $m$, then $ x= \phi_t(y)$ for some $y \in F_m \subset  S^+(\dM)\setminus \trap^+_\Omega(\dM)$
and some $t \in [0, t(x)]$. Then $y \notin \trap^+_\Omega(\dM)$ gives $t(y) < \infty$ and so $t(x) < \infty$, too. Also, if $x = (q,v)$,
then $t(q,-v) = t < \infty$. Thus, $x\notin \trap(\Omega)$. This proves that $G_m \subset S(\Omega) \setminus \trap (\Omega)$ for 
all $m$. Next, let $x = (q,v) \in S(\Omega) \setminus \trap (\Omega)$. Then $t = t(q,-v) < \infty$, so $x = \phi_t(y)$ for some
$y \in S^+(\dM)$. Moreover, $t(y) = t(x) + t < \infty$, so $y \in S^+(\dM)\setminus \trap^+_\Omega(\dM)$, and by (2.7) we have
$y \in F_m$ for some $m \geq 1$. 
Moreover $0 \leq t \leq t(y)$. Thus, $x \in G_m$. This proves (2.8).

Notice that all sets $G_m$ are closed in $S(\Omega)$. Indeed, given $m$, a standard argument (see e.g. the proof of Lemma 3.1 below) 
shows that the function $t(x)$ is uniformly bounded on the compact set $F_m$. Let $\{x_p\}_{p=1}^\infty$ be a sequence in $G_m$
with $x_p \to x$ as $p\to\infty$. Then $x_p = \phi_{t_p}(y_p)$ for some $t_p \in [0,T]$ and some $y_p \in F_m$, where
$T = \sup_{y\in F_m} t(y)$. Taking an appropriate sub-sequence, we may assume $y_p \to y$ and $t_p \to t$ as $p \to \infty$.
Since $F_m$ is compact, $y \in F_m$. Now a simple continuity argument gives $x = \phi_t(y)$, so $x \in G_m$. 

Thus, $G_m$ is closed in $S(\Omega)$ and so $V_m = S(\Omega) \setminus G_m$ is open in $S(\Omega)$. Moreover, it
is clear that $\phi_t(x) \in V_m$ for all $x\in V_m$ and all  $t \in [0,t(x)]$. Let $\chi_{G_m}$ be the {\it characteristic function} of $G_m$ in $S(\Omega)$.
Applying Lemma 2.1 to $V = V_m$ and $f$ replaced by $f_m = f\cdot \chi_{G_m}$,  we get
$$\int_{S(\Omega)\setminus \trapf(\Omega)} f_m(x) \; d\lambda(x)  = \int_{S^+(\partial M)} \left( \int_0^{t(x)} f_m(\phi_t(x)) \;dt \right)\, d\mu(x) .$$
Extending $f_m$ as $0$ on $V_m$, we get a sequence of measurable functions on $S(\Omega)$ with $|f_m| \leq |f|$ for all
$m$ and $f_m(x) \to f(x)$ for all $x\in S(\Omega)$. Since $|f|$ is integrable,  Lebesgue's Theorem now implies
$$\int_{S(\Omega)\setminus \trapf(\Omega)} f (x) \; d\lambda(x)  = \int_{S^+(\partial M)} \left( \int_0^{t(x)} f (\phi_t(x)) \;dt \right)\, d\mu(x) ,$$
which proves the theorem.
\endofproof


\section{Proof of Theorem 1.2}\label{sec3}
\renewcommand{\theequation}{\arabic{section}.\arabic{equation}}
\setcounter{equation}{0}

 It is enough to prove the more general part (b).

Let $M$ and $\Omega$ satisfy the assumptions of Theorem 1.2, so (A) holds. Let $\ep > 0$.
Assume that there is no open neighbourhood $U$ of $\id$ in $C^k(\dk,M)$  such that for 
every $F \in U$ we have $|\lambda (\trap(\Omega_m)) - \lambda(\trap(\Omega) | < 2\ep$. Then there exists a
sequence $\{F_m\}_{m=1}^\infty \subset C^k(\dk, M)$ with $F_m \to \id$ as $m\to \infty$ in the $C^k$ Whitney topology
such that 
\begin{equation}
|\lambda (\trap(\Omega_m)) - \lambda(\trap(\Omega) | \geq 2 \ep
\end{equation}
for all $m$, where we set $\Omega_m = \Omega_{F_m}$ for brevity. 
Set $K_m = \overline{M \setminus \Omega_m}$; then $\dk_m = F_m(\dk)$. Since all obstacles $K_m$ are in the interior
of $M$, we have $\Omega_m \subset M$, so the Liouville measure $\lambda = dqdv$ is well-defined on $S(\Omega_m)$. 
The function $t(x)$ and the billiard trajectory $\gamma^+(x)$ for $\Omega_m$ will be denoted by $t_m(x)$ and $\gamma^+_m(x)$, 
respectively. Set
$$\tt' = \trap^+_\Omega(\dM) \cup \cup_{m=1}^\infty \trap^+_{\Omega_m} (\dM) .$$
Then $\mu(\tt') = 0$ by (1.2). Theorem 1.1 with $f = 1$ implies
$$\lambda(S(\Omega)\setminus \trap(\Omega)) = \int_{S^+(\dM)} t(x)\, d\mu(x) =  \int_{S^+(\dM)\setminus \tt'} t(x)\, d\mu(x) ,$$
and
$$\lambda(S(\Omega_m)\setminus \trap(\Omega_m))  = \int_{S^+(\dM)} t_m(x)\, d\mu(x) =  \int_{S^+(\dM)\setminus \tt'} t_m(x)\, d\mu(x) .$$
Since $\lambda(S(\Omega_m)) \to \lambda(S(\Omega))$ as $m \to \infty$, we may assume that 
$|\lambda(S(\Omega_m)) - \lambda(S(\Omega))| < \ep$ for  all $m \geq 1$. Combining this with (3.1) and
the above  equalities, we get
\begin{equation}
\left|  I_m -  I \right| \geq \ep
\end{equation}
for all $m \geq 1$, where we set for brevity
$$I_m = \int_{S^+(\dM)\setminus \tt'} t_m(x)\, d\mu(x) \quad , \quad I = \int_{S^+(\dM)\setminus \tt'} t (x)\, d\mu(x) .$$

For every $m \geq 1$ denote by $\tt''_m$ the set of those $x\in S^+(\dM)$ such that $\gamma^+_m(x)$ has a tangent point to $\do_m$, and let
$\tt''_0$ be the set of those $x\in S^+(\dM)$ such that $\gamma^+(x)$ has a tangent point to $\do$.
As remarked earlier, we have $\mu(\tt''_m) = 0$ for all $m \geq 0$. Thus, for the set
$$\tt'' = \cup_{m=0}^\infty \tt''_m$$
we have $\mu(\tt'') = 0$, and therefore $\mu(\tt' \cup \tt'') = 0$. Let
$$D_1 \subset D_2 \subset \ldots \subset D_r \subset \ldots \subset S^+(\dM)\setminus (\tt'\cup \tt'')$$
be compact sets so that $\mu(D_r) \nearrow \mu(S^+(\dM))$ as $r \to \infty$. 
Set 
$$\Ir_m = \int_{D_r} t_m(x)\, d\mu(x) \quad , \quad \Ir = \int_{D_r} t (x)\, d\mu(x) $$
for brevity. Then 
\begin{equation}
\Ir_m \nearrow I_m \quad, \quad \Ir \nearrow I
\end{equation}
as $r \to \infty$, for any fixed $m$ in the first limit. It follows from this that there exists an integer $r_0$ so that
\begin{equation}
I -  \frac{\ep}{3} < \Ir \leq I
\end{equation}
for all $r \geq r_0$.

Next, consider an arbitrary point $x = (q,v) \in S^+(\dM) \setminus (\tt' \cup \tt'')$. Then $x\notin \trap^+_\Omega(\dM)$ gives $t(x) < \infty$,
while $x\notin \tt''_0$ shows that $\gamma^+(x)$ is a simply (transversally) reflecting trajectory in $\Omega$ with no tangencies to $\do$. Thus,
{\it the number $j_0(x)$ of reflection points} of $\gamma^+(x)$ is finite, and $\gamma^+(x)$ has no other common points with $\do$. 

\bs

\noindent
{\bf Lemma 3.1.} {\it Let $r \geq 1$. There exists an integer $m_0 = m_0(r) \geq 1$ such that for every $x \in D_r$ and every integer $m \geq m_0$
the trajectory $\gamma^+_m(x)$ has at most $j_0(x)$ common points\footnote{Clearly all of them must be simple reflection points,
since $x \notin \tt'_m \cup \tt''_m$.} with $\do_m$.}

\bs

\noindent
{\bf Proof of Lemma 3.1.} Fix $r \geq 1$ and assume there exist arbitrarily large $m$ such that $\gamma^+_m(x_m)$ has at least $j_0(x_m) +1$ reflection points
for some $x_m \in D_r$. Choosing a subsequence, we may assume that the latter is true for all $m \geq 1$ and that $x_m \to x \in D_r$ as $m \to\infty$.
Since both $\gamma^+(x_m)$ and $\gamma^+(x)$ are simply  reflecting, it follows that $j_0(x_m) = j_0(x)$ for all sufficiently large $m$; we will assume
this is true for all $m$. Set $j_0 = j_0(x)$. Thus, for every $m$ the trajectory $\gamma^+_m(x_m)$ has at least $j_0 +1$ reflection points. 
Let $q_1, \ldots, q_{j_0} \in \dk$ be the successive reflection points of $\gamma^+(x)$. Set $q_{j_0+1} = g_1$ for convenience. Notice that some of the points $q_j$
may coincide, however $q_{j+1} \neq q_j$ for all $j$. For each $j = 1, \ldots, j_0$ choose a small open neighbourhood $V_j$ of $q_j$ on $\dk$ so that 
$\overline{V_j} \cap \overline{V_{j+1}} = \e$ for all $j = 1,\ldots,j_0$.  Since $F_m \to \id$ in the $C^k$ Whitney topology, a simple
continuity argument shows that for sufficiently large $m$ we have  $\overline{F_m(V_j)} \cap \overline{F_m(V_{j+1})} = \e$ for all $j = 1, \ldots,j_0$. Assuming that
the neighbourhoods $V_j$ are chosen sufficiently small and $m$ is sufficiently large, it follows that for each $j = 1, \ldots,j_0$, $\gamma^+_m(x_m)$ has an
unique reflection point $q^{(m)}_j$ in $F_m(V_j)$. Let $q^{(m)}_{j_0+1}$ be a reflection point of $\gamma^+_m(x_m)$ which is not in $\cup_{j=1}^{j_0} F_m(V_j)$;
such a point exists, since by assumption the trajectory $\gamma^+_m(x_m)$ has at least $j_0 +1$ reflection points. Choosing an appropriate sub-sequence again,
we may assume that $q^{(m)}_{j_0+1} \to q_{j_0+1} \in \dk$ as $m \to \infty$. It is now clear that $q_{j_0+1}$ is a common point of $\gamma^+(x)$ and $\dk$
which does not belong to $\cup_{j=1}^{j_0} V_j$, so $q_{j_0+1} \neq q_j$ for all $j = 1, \ldots, j_0$. This is a contradiction which proves the lemma.
\endofproof

\ms

We now continues with the proof of Theorem 1.2.

According to (3.2), we either have $I_m \leq I - \ep$ for infinitely many $m$, or
$I_m \geq I+\ep$ for infinitely many $m$.

\ms

\noindent
{\bf Case 1.} $I_m \geq I + \ep$ for infinitely many $m$. Considering an appropriate subsequence, we may assume $I_m \geq I + \ep$ for 
all  $m \geq 1$.

Let  $r \geq r_0$ so that (3.4) holds for $r$, and let $m_0 = m_0(r)$ be as in Lemma 3.1.
Then the compactness of $D_r$ and a simple  continuity argument show that 
$$i_0 = \sup_{x\in D_r} j_0(x) < \infty .$$
This and Lemma 3.1 now imply 
$$t_m(x) \leq D\, i_0 \quad, \quad m \geq m_0 \;, \; x\in D_r .$$
Another simple continuity argument shows that for $x\in D_r$ the only possible limit point of the sequence $\{ t_m(x)\}_{m=1}^\infty $ is $t(x)$, so we must have
$\lim_{m\to\infty} t_m(x) = t(x)$ for all $x\in D_r$. This is true for all $r$, so $\lim_{m\to\infty} t_m(x) = t(x)$ for all $x\in S^+(\dM) \setminus (\tt' \cup\tt'')$.
Setting
$$\ttt_m(x) = \sup_{m' \geq m} t_{m'}(x) \quad, \quad x\in  S^+(\dM) \setminus (\tt' \cup\tt'') ,$$
we get a sequence of measurable functions with  $\ttt_m(x) \searrow t(x)$ on 
$S^+(\dM) \setminus  (\tt' \cup\tt'')$.
By Lebesgue's Theorem (or Fatou's Lemma),
\begin{eqnarray*}
\lim_{m\to\infty} \int_{S^+(\dM)} \ttt_m(x) \, d\mu(x) 
& = & \lim_{m\to\infty} \int_{S^+(\dM)\setminus (\tt'\cup \tt'')} \ttt_m(x) \, d\mu(x)\\
& = & \int_{S^+(\dM)\setminus (\tt'\cup \tt'')} t(x) \, d\mu(x) =    \int_{S^+(\dM)} t(x) \, d\mu(x).
\end{eqnarray*}
Thus, there exists $m_1 \geq m_0$ such that
$$\int_{S^+(\dM)} \ttt_m(x) \, d\mu(x)  <  \int_{S^+(\dM)} t(x) \, d\mu(x) + \frac{\ep}{3} = I + \frac{\ep}{3} $$
for all $m \geq m_1$. Since $t_m(x) \leq \ttt_m(x)$, it follows that
\begin{equation}
I_m = \int_{S^+(\dM)} t_m(x) \, d\mu(x)  \leq \int_{S^+(\dM)} \ttt_m(x) \, d\mu(x)<  I + \frac{\ep}{3}
\end{equation}
for $m \geq m_1$, which is a contradiction with our assumption that $I_m \geq I + \ep$ for  all  $m \geq 1$.

\ms

\noindent
{\bf Case 2.} $I_m \leq I - \ep$ for infinitely many $m$. Considering an appropriate subsequence, we may assume $I_m \leq I - \ep$ for 
all  $m \geq 1$. Combining this with (3.3) and (3.4) gives
\begin{equation}
\Ir_m \leq I_m \leq I- \ep < \Ir + \frac{\ep}{3} - \ep = \Ir - \frac{2\ep}{3}
\end{equation}
for all $r \geq r_0$ and all $m \geq 1$.

Next, {\bf fix an arbitrary integer $r \geq r_0$} so that (3.4) holds for $r$, and let $m_0 = m_0(r)$ be as in Lemma 3.1. 
As in Case 1,
the compactness of $D_r$ and a simple  continuity argument show that 
$$i_0 = \sup_{x\in D_r} j_0(x) < \infty ,$$
while Lemma 3.1 implies 
$$t_m(x) \leq D\, i_0 \quad, \quad m \geq m_0 \;, \; x\in D_r .$$
As in Case 1, for any $x\in D_r$ we must have $\lim_{m\to\infty} t_m(x) = t(x)$ for all $x\in D_r$. 

Moreover, for the fixed $r$, this convergence is uniform. Indeed, if this is not true, then there exist $\delta > 0$
and infinite sequences $\{ x_s\} \subset D_r$ and $1 \leq m_1 < m_2 < \ldots < m_s < \ldots$ such that $|t_{m_s} (x_s) - t(x_s)| \geq \delta$ for all $s$.
Using the compactness of $D_r$, we may assume that $x_s \to x\in D_r$ as $s \to \infty$. Also, since $\{t_{m_s}(x_s)\}$ is a bounded sequence, 
we may assume that $t_{m_s}(x_s) \to t\in [0, D\, i_0]$ as $s \to \infty$. 
Since the trajectory $\gamma^+_{m_s}(x_s)$ has at most $i_0$ reflection points, choosing an appropriate subsequence again, we may assume that the
billiard trajectory  $\gamma^+_{m_s}(x_s)$ has the same number $p$ of reflections points for all $s \geq 1$. 
Let $y_1(s), y_2(s), \ldots, y_p(s)$ be the successive reflection points of $\gamma^+_{m_s}(x_s)$. By compactness, choosing appropriate subsequences again,
we may assume that $y_i(s) \to y_i$ as $s \to \infty$ for all $i = 1, \ldots, p$. It is now clear that $y_1, \ldots, y_p \in \dk$ and these are the successive reflection
points of a billiard trajectory in $\Omega$. Since $x_s \to x$ as $x \to \infty$, we must have that $y_1, \ldots, y_p$ are the successive reflection points of
$\gamma^+(x)$. In particular, we must have $t_{m_s}(x_s) \to t(x)$ as $s \to \infty$. However, $t(x_s) \to t(x)$ as  well, so it follows that
$|t_{m_s} (x_s) - t(x_s)| < \delta$ for all sufficiently large $s$; a contradiction with our assumption.

Thus, $t_m(x) \to t(x)$ as $m \to \infty$ uniformly for $x\in D_r$. This implies that there exists $m_1 \geq m_0$ such that
$$|\Ir_m - I| = \left| \int_{D_r} t_m(x) \, d\mu(x) - \int_{D_r} t(x)\, d\mu(x) \right| < \frac{\ep}{3} $$
for all $m \geq m_1$. In particular we have $\Ir_m > I - \frac{\ep}{3}$ for all $m \geq m_1$, which is a contradiction 
with (3.6). 

\ms

In this way we have show that (3.2) cannot hold for (infinitely many)  $m \geq 1$.
This completes the proof of the theorem.
\endofproof

\section{Scattering by obstacles in $\R^n$}\label{sec3}
\renewcommand{\theequation}{\arabic{section}.\arabic{equation}}

\setcounter{equation}{0}

Here we consider the case when $\tM = \R^n$ with the standard Riemannian metric for some $n \geq 2$, and
$K$ is a compact subset of ${\R}^n$ whose boundary $\partial K$ is a $C^k$ manifold of dimension $n-1$ for
some $k \geq 2$.  
We assume that ${\R}^n\setminus K$ is connected. Let $M$ be a large closed ball in $\R^n$ containing $K$ 
in its interior.  As in Sect. 1,  
$$\Omega = \overline{M\setminus K}$$ 
has a smooth boundary $\do = S_0 \cup \dk$,  where $S_0$ is the {\it boundary sphere} of $M$.

In the present case the scattering rays from Sect. 1 are simply billiard trajectories in the exterior of $K$ that 
come from infinity, enter $M$ at some point $q\in S_0$ with a certain direction $v \in \sn$ and after a 
time $t(q,v)$ spent in $\Omega$, 
leave $M$ and go to infinity. Then $t(q,v)$ is what we called the travelling time of $x = (q,v)$ in Sect. 1.

It is a natural problem to try to recover information about the obstacle $K$ from measurements along billiard 
trajectories (generalised geodesics) in the exterior of $K$. As we mentioned in the Introduction, 
problems of this kind have been considered for a rather long time in Riemannian geometry and more
recently in scattering by obstacles, as well. Reconstructing $K$ in practice from the travelling times data appears 
to be a rather difficult  problem, although in relatively simple cases there is enough scope to achieve this -- see for example Sect. 4 in 
\cite{NS2} which describes how to recover a planar obstacle $K$ which is a disjoint union of two 
strictly convex domains.

In the case considered in this section, condition (A) from Sect. 1 is always satisfied. 
A point 
$x = (q,v)\in  S(\Omega)$ is {\it trapped} if either  its forward billiard trajectory $\gamma^+(x)$ or 
its backward trajectory $\gamma^+(q, -v)$ is infinitely long. That is, either the billiard 
trajectory in the exterior of $K$ issued from $q$ in the direction
of $v$ is bounded (contained entirely in the ball $M$) or the one issued from $q$ in the direction of $-v$ is bounded.


What concerns the problem of obtaining information about $K$ from travelling times $t(x)$, 
Theorem 1.1 provides some general information. 
In particular when $f = 1$ we get the following consequence.

\bs

\noindent
{\bf Theorem 4.1.} {\it Assume that  the set $\trapk$ of trapped points in $\Omega_K$ has Lebesgue measure zero. 
Then 
\begin{equation}
\Vol_n(K) = \Vol_n(M)  - \frac{1}{\Vol_{n-1}(\sn)}\; \int_{S^+(S_0)} t(x) \; d\mu (x) ,
\end{equation}
where $\Vol_n(K)$ is the standard volume of $K$ in $\R^n$ and $\Vol_{n-1}(\sn)$ 
is the standard  $(n-1)$-dimensional volume (surface area) of $\sn$}.

\bs

\noindent
{\bf Proof.} It follows from Theorem 1.1 with $f = 1$ that
$$\Vol(S(\Omega))  =  \int_{S^+(S_0)} t(x) \; d\mu (x) .$$
Combining this with  
$\Vol(S(\Omega)) = \Vol_n(\Omega) \, \Vol_{n-1}(\sn)$ and  
$\Vol_n(K) = \Vol_n(M) - \Vol_n(\Omega)$ proves (4.1).
\endofproof

\bs


\noindent{\bf Remarks.}
(a) Formula (4.1) shows that when $\trapk$ has measure zero, from travelling times data we can recover
the volume of $K$. That is, without seeing $K$ and without any preliminary
information about $K$, just measuring travelling times
of a certain kind of signals incoming through points on the sphere $S_0$ and outgoing through
points on $S_0$, we can compute the amount of mass in $K$, i.e. the volume of $K$. Apart from that,
it appears that (4.1) could be useful in numerical approximations of the volume of $K$.

\ms

(b) In Theorem 4.1 we only used the trivial function $f = 1$. Naturally, one would expect that
using Theorem 1.1 for a large family of functions $f$  would bring much more significant information
about the obstacle $K$. It is already known (see e.g. \cite{St3} and \cite{NS2}) that a
solid amount of information about $K$ is recoverable from travelling times. However by
means of formula (1.3) it might be possible to get such information in a more explicit way.

\bs

In simple cases when $K$ is a disjoint union of connected pieces of roughly the same size and shape, 
we can estimate the number $k$ of these pieces. Naturally, $k$ can be a very large number\footnote{E.g. consider 
a model that resembles the molecules in a gas container.}.

\bs

\noindent
{\bf Example.} Assume that $K$ is a disjoint union of $k$ balls of the same radius $a > 0$, where $k \geq 1$ 
is arbitrary (possibly a large number). Suppose that we know $a$ from some preliminary information. 
Then measuring travelling times $t(x)$ for a relatively large number of points $x = (q,v) \in S(S_0)$ we get 
an approximation of the integral  $ \int_{S^+(S_0)} t(x) \; d\mu (x)$, and therefore by means of (4.1), we obtain 
an approximate value for the number $k$ of connected components of $K$. The precise formula (assuming we can 
measure almost all travelling times) follows from (4.1):
$$ k = \frac{\Vol_n(K)}{\pi^{n/2}  a^n/\Gamma(n/2+1)} =  \frac{R^n}{a^n}
 - \frac{\Gamma(n/2)\; \Gamma(n/2+1)}{2 \pi^n a^n}\; \int_{S^+(S_0)} t(x) \; d\mu (x)  ,$$
where $R$ is the radius of the large ball $M$ containing $K$, and $\Gamma$ is Euler's Gamma function.

\section{A Corollary}\label{sec5}
\renewcommand{\theequation}{\arabic{section}.\arabic{equation}}

\setcounter{equation}{0}

Given an integer $k \geq 0$, recall from Sect. 2 that (in the case $V = \e$ in the proof of Lemma 2.1) 
$\Gamma_k$ is {\it the set of all $x\in S^+(\partial M)$  
such that $t(x) < \infty$, $\gamma^+(x)$ contains no gliding segments on the boundary $\do$ and 
has exactly $k$ reflections at} $\dk$. As another consequence of Theorem 1.1 we get the
following, the first part of which concerns the general situation considered in Sects. 1 and 2, while
the second deals with a special kind of obstacles in $\R^n$.

\bs


\noindent
{\bf Corollary 5.1.} {\it Let $D = \diam(M)$.}

\ms

(a) {\it Under the assumptions of Theorem {\rm 1.1} we have}
$$\frac{1}{D}  \Vol(S(\Omega)) \leq \sum_{k=0}^\infty (k+1) \, \mu(\Gamma_k) .$$

\ms

(b) {\it Let $K$ be a finite disjoint union of strictly convex domains in $\R^n$ ($n \geq 2$)
with smooth boundaries, and let $d > 0$ be a constant such that the minimal distance between distinct connected components 
of $K$ is not less than $d$ and $d\leq \dist(K, \partial M)$, where $M$ is a ball of diameter $D$ containing $K$ in its interior. Then
$$\frac{1}{D} \Vol(S(\Omega)) \leq \sum_{k=0}^\infty (k+1) \, \mu(\Gamma_k) \leq 
\frac{1}{d} \Vol(S(\Omega)) .$$
In particular, there exists a constant $C > 0$ such that
$$\di \mu(\Gamma_k) \leq \frac{C}{k+1} $$
for all $k \geq 0$.}

\bs

\noindent
{\bf Proofs.} (a) By Theorem 1.1 with $f = 1$ we get
\begin{eqnarray*}
\Vol(S(\Omega)) 
& = &  \int_{S^+(S_0)} t(x) \; d\mu (x) = \sum_{k=0}^\infty \int_{\Gamma_k} t(x)\; d\mu (x)
\leq \sum_{k=0}^\infty \int_{\Gamma_k} (k+1) D\; d\mu (x)\\
& = & D \; \sum_{k=0}^\infty (k+1) \, \mu(\Gamma_k) ,
\end{eqnarray*}
which proves part (a).

\ms

(b) First, notice that $\tM = \R^n$, $M$, $K$ and $\Omega = \overline{M\setminus K}$
satisfy the condition (A) from Sect. 1. Morever, it follows from Proposition 2.3 in \cite{St1}
and Proposition 5.1 in \cite{St2} that the set Trap$(K)$ of trapped points has Lebesgue
measure zero. By the nature of $K$ and $M$, no billiard trajectory in $S(\Omega)$
contains non-trivial gliding segments on $\do$. Thus, for any non-trapped point $x\in S(\Omega)$
the trajectory $\gamma^+(x)$ has only finitely many reflection points.

Hence Theorem 1.1 is applicable. Using it again with $f = 1$ as in the proof of part (a),
this times estimating from below $t(x) \geq (k+1) d$, proves the assertion.
\endofproof

\end{document}